\documentclass[reqno]{amsart} 
\begin{document} 

\title[\hfilneg \hfil Sturm-Liouville problems with prime spectrum]
{A note on Sturm-Liouville problems whose spectrum is the set of prime numbers} 

\author[ A. B. Mingarelli]
{Angelo B. Mingarelli}  % in alphabetical order

\address{School of Mathematics and Statistics\\ 
Carleton University, Ottawa, Ontario, Canada, K1S\, 5B6}
\email[A. B. Mingarelli]{amingare@math.carleton.ca}

\date{}
\thanks{Submitted August ??, 2011.}
\thanks{The author is partially supported by a grant from NSERC Canada}

\begin{abstract}
 We show that there is no classical regular Sturm-Liouville problem on a finite interval whose spectrum 
consists of infinitely many distinct primes numbers. In particular, this answers in the negative a question raised by Zettl in his book, \cite{zet}. We also show that there {\it may} exist such a problem if the parameter dependence is nonlinear.
\end{abstract}

\maketitle
\numberwithin{equation}{section}
\newtheorem{theorem}{Theorem}[section]
\newtheorem{corollary}[theorem]{Corollary}
\newtheorem{example}[theorem]{Example}

\section{Introduction}
In this note a Sturm-Liouville equation on a finite interval $[a,b]$ is defined by a second order real differential expression of the form
\begin{equation}\label{1.1}
-(p(x)y^{\prime})^{\prime} + q(x)y = \lambda r(x) y
\end{equation}
where $p, q, r : [a,b] \to \mathbb{R}$ and $1/p, q, r \in L[a,b]$, and $\lambda$ is a generally complex parameter. By a {\it solution} of \eqref{1.1} we will mean, as is customary, a function $y$ defined and absolutely continuous on $[a,b] = I$ such that $(py^{\prime})(x)$ is also absolutely continuous and $y(x)$ satisfies the differential equation \eqref{1.1} almost everywhere on $I$. In this setting one can allow $p(x)$ to have infinite values but only on a set of measure zero.

The classical regular Sturm-Liouville problem (SLP) associated with \eqref{1.1} consists in finding those values of $\lambda \in \mathbb{C}$ such that \eqref{1.1} has a non-trivial solution satisfying the {\it separated homogeneous boundary conditions}
\begin{eqnarray}
\label{1.2}
y(a)\cos \alpha - (py^{\prime})(a)\sin \alpha =0,\\
\label{1.3}
y(b)\cos \beta - (py^{\prime})(b)\sin \beta =0,
\end{eqnarray}
where $\alpha \in [0,\pi)$ and $\beta \in (0,\pi]$. Of course, this problem has a very long history dating back to Sturm's original contributions in the 1830's. It has been developed in many different directions but we shall not delve on this at the moment.

Still greater generality can be obtained by allowing $p(x)$ to be identically infinite on subintervals of $I$. In this case one needs to rewrite \eqref{1.1} as a vector system in two dimensions, e.g., 
\begin{eqnarray}
\label{1.4}
u^{\prime} = - v/p, \\
\label{1.5}
v^{\prime} = (\lambda r - q) u.
\end{eqnarray}

This defines a problem of {\it Atkinson-type}. The boundary conditions \eqref{1.2}-\eqref{1.3} now take the form 
\begin{eqnarray}
\label{1.6}
u(a)\cos \alpha + v(a)\sin \alpha =0,\\
\label{1.7}
u(b)\cos \beta + v(b)\sin \beta =0,
\end{eqnarray}
where $\alpha \in [0,\pi)$ and $\beta \in (0,\pi]$.

The advantage of using the formulation \eqref{1.4}-\eqref{1.5} is that it can be used to study three-term recurrence relations as well, see [\cite{bvp}, Chapter 8] and \cite{abm} for more details. Here, the lack of assumptions on the signs of the coefficients allows for the most generality. However, in this very case things can get pretty bad since there are coefficients and a corresponding interval $I$ with the property that the Dirichlet problem for \eqref{1.1} on $I$ has discrete spectrum filling the whole complex plane (see \cite{atm}, \cite{abm4}). For more information on Sturm-Liouville problems of Atkinson-type see \cite{bvp} and \cite{zet}.

The aim of this note is to answer, in part, Problem IV raised by Zettl in [\cite{zet}, p. 299]. We restate the problem here for ease of reference:

\begin{center}``IV: {\it Find a SLP whose spectrum is the primes.}"\end{center}

It is then pointed out that given any {\it finite} set of distinct real numbers, a SLP of Atkinson-type can be found whose spectrum is precisely that set, [\cite{zet}, p. 299].

On the basis of spectral asymptotics and the Prime Number Theorem we show that there does not exist a regular Sturm-Liouville problem whose spectrum consists of infinitely many prime numbers, thus answering said Problem IV in the negative. However, we show that there {\it may} exist such a problem if the parameter dependence in \eqref{1.1} (or \eqref{1.2}-\eqref{1.3}) is {\it nonlinear}. Nevertheless, it is not clear how to choose the coefficients so that the primes are generated even in such nonlinear cases.

\section{The main results and discussions}

In the sequel the set of rational primes refers to the usual set of prime numbers in the rational number field. For basic results about primes and their distribution we refer the reader to Hardy and Wright \cite{hw}. The following two results are, for the most part,  independent of sign conditions on the coefficients in either \eqref{1.1} or \eqref{1.4}-\eqref{1.5} (see \cite{atm} for specific details).
\begin{theorem}
\label{th1}
There does not exist a regular SLP \eqref{1.1} with separated boundary conditions \eqref{1.2}-\eqref{1.3} whose spectrum is an infinite set of distinct rational primes.
\end{theorem}

Indeed, more is true. Theorem~\ref{th1} has a counterpart for SLP of Atkinson-type. \begin{theorem}
\label{th2}
There does not exist a regular SLP of Atkinson-type \eqref{1.4}-\eqref{1.5} with separated boundary conditions \eqref{1.6}-\eqref{1.7} whose spectrum is an infinite set of distinct rational primes.
\end{theorem}
It follows that Zettl's Problem IV as stated is unsolvable within this framework. 

Now consider the Dirichlet problem associated with \eqref{1.1} on $[0,1]$ with nonlinear dependence 
\begin{equation}\label{1.8}
-y^{\prime\prime} + q(x)y = \left (\frac{\pi \lambda}{\log \lambda}\right )^2\, y.
\end{equation}
The basic asymptotic estimate for the rational primes states that if $p_n$ denotes the $n$-th prime, then 
\begin{equation}\label{2.3}p_n \sim n\, \log n, \quad n\to \infty,\end{equation}
the latter being a consequence of the Prime Number Theorem, cf., [\cite{hw}, p.10 and p.346].

Now the positive eigenvalues of the Dirichlet problem for \eqref{1.8}  (whose existence is essentially a consequence of the Atkinson-Mingarelli Theorem, \cite{atm}) have the {\it correct} asymptotics, i.e., 
\begin{equation} \lambda_n \sim n\, \log n, \quad n\to \infty.
\end{equation}

Choosing the simplest case (i.e., $q(x)=0$) we get that the eigenvalues $\lambda_ n$, for $n \geq 1$, of the Dirichlet problem on $[0,1]$ of \eqref{1.8} must be the roots of the equation
$$\frac{\lambda}{\log \lambda} = n,$$ 
from which we deduce the approximate values of $\lambda_n$, for large $n$, namely:
$$\lambda_n = n\log n + n\log \log n + n \log\log\log n + \cdots$$
Even so, without any specific reference to prime numbers, the first two terms on the right of the preceding equation already agree with Ces\`{a}ro's  approximate value for the $n$-th prime \cite{ec}, i.e.,
$$p_n = n\log n + n\log \log n - n + n\frac{(\log\log n - 2)}{\log n} - \cdots$$

It is thus {\it conceivable}, but by no means obvious, that by modifying $q(x) \in L[0,1]$ on infinitely many subintervals on each of which $q(x)$ takes on constant values, we may come up with an asymptotic expansion for $p_n$ for large $n$.   If so, we can hope to find a suitable $q$ such that the Dirichlet spectrum of \eqref{1.8} agrees with the set of all rational primes. We emphasize that the existence of such a function $q(x)$ is an open question.

\section{Proofs}
\begin{proof} (Theorem~\ref{th1}) We assume that $p(x)$ has a finite number of turning points and  that $(r(x)/p(x))_+ $ is not a.e. equal to zero. Here, $(\cdot)_+$   denotes the positive part of the function in question. Then, by a direct application of the Atkinson-Mingarelli Theorem [\cite{atm}, Theorem~2.3], we get that the positive eigenvalues, $\lambda_n^+$, of \eqref{1.1}-\eqref{1.3} satisfy
\begin{equation}
\label{3.1}
\lambda_n^+ \sim \dfrac{n^2\,\pi^2}{\int_{a}^{b} \sqrt{\left ({r(x)}/{p(x)} \right)_+ \, dx}}.
\end{equation}
It is now clear that \eqref{2.3} and \eqref{3.1} are incompatible if the primes (or an infinite subsequence of distinct such) were in fact the eigenvalues of such a problem.
\end{proof}

\begin{proof} (Theorem~\ref{th2})  Although the Atkinson-Mingarelli Theorem does not extend immediately to SLP of Atkinson-type one can resort to a simple order of magnitude argument. By assumption we take it that the spectrum is an infinite set at the outset (more on this below). We argue as in [\cite{bvp}, pp.206-207] with minor changes. Assuming that for the moment $u(a), v(a)$ have fixed initial values, it is a known fact that the resulting unique solution, $u(x,\lambda)$, and indeed $v(x, \lambda)$, are analytic for all complex $\lambda$ and thus each is an entire function of $\lambda$ (see also \cite{abm2} and \cite{abm3} for extensions of these results). 

Now, following  [\cite{bvp}, pp.206-207] a straightforward calculation shows that, even for indefinite coefficients, we still have
$$ \bigg | \frac{d}{dx} \log \{ |\lambda| |u|^2 + |v|^2\} \bigg | \leq \sqrt{|\lambda|}(|r| + |1/p|) + |q|/ \sqrt{|\lambda|}.$$
Integrating the latter over $[a,b]$, using the integrability assumptions on the coefficients, and exponentiating, we obtain the estimates (see [\cite{bvp}; p. 206, eq. (8.2.5)])
$$ u(x,\lambda), v(x,\lambda) = \Large{O}\{exp(const. \sqrt{|\lambda|)}\},$$
thus showing that, as entire functions of $\lambda$, each must be of order not exceeding $1/2$. It follows from the theory of entire functions that the infinite sequence of complex zeros, $\lambda_n$, of \eqref{1.7} (that define the eigenvalues) grow at such a rate that, for any $\varepsilon > 0$,
$$\sum_{n=1}^\infty \frac{1}{|\lambda_n|^{\frac{1}{2} + \varepsilon}} < \infty.$$
However, for sufficiently small $\varepsilon > 0$, \eqref{2.3} implies that the series of primes (or any infinite subsequence of distinct primes) must satisfy 
$$\sum_{n=1}^\infty \frac{1}{|p_n|^{\frac{1}{2} + \varepsilon}} = \infty.$$
\end{proof}
This contradiction completes the proof.

\end{document}